%

\magnification=1200
\input amstex
\documentstyle{amsppt}

\define\be{\overline{\eta_\alpha}}
\define\eb{\bar\eta}
\define\res{\restriction}
\define\supp{\text{supp}}
\TagsAsMath

\topmatter
\title
Subgroups of the Baer-Specker Group\\
with Few Endomorphisms but Large Dual
\endtitle
\rightheadtext{Subgroups of the Baer-Specker Group}
\author
Andreas Blass and R\"udiger G\"obel
\endauthor
\address
Mathematics Dept., University of Michigan, Ann Arbor, 
MI 48109, U.S.A.
\endaddress
\email
ablass\@umich.edu
\endemail
\address
Fachbereich 6, Mathematik, Universit\"at Essen,
Postfach, 45117 Essen, Germany
\endaddress
\email
 math100\@vm.hrz.uni-essen.de
\endemail
\thanks Partially supported by NSF grant DMS-9204276 and NATO grant
LG921395.
\endthanks
\abstract
Assuming the continuum hypothesis, we construct a pure subgroup $G$ of
the Baer-Specker group $\Bbb Z^{\aleph_0}$ with the following
properties.  Every endomorphism of $G$ differs from a scalar
multiplication by an endomorphism of finite rank.  Yet $G$ has
uncountably many homomorphisms to $\Bbb Z$.
\endabstract
\endtopmatter
\document

\head
{\S 1 Introduction}
\endhead

The Baer-Specker group, which is the additive group $P =\Bbb Z^{\aleph
_0}$ of all integer-valued sequences, and its obvious generalization
$\Bbb Z^\kappa $, the cartesian product of $\Bbb Z$ for any infinite
cardinal $\kappa $, have attracted continuous research over many years.
This is not too surprising because many questions on abelian groups,
like the famous Whitehead problem, can be stated in terms of subgroups
of such products. Recall that

\noindent
$(\dag )$\quad {\sl $P$ is $\aleph_1$-free,
i.e. all countable subgroups of $P$ are free}

\noindent
by a result of Baer, cf. Fuchs [15, Vol.~1, p.~94, Theorem~19.2].
Besides the Whitehead problem, many other properties of these products
turn out to depend on the set theoretic assumptions, which make these
groups a favored playground for logicians as well as algebraists.
Examples of results of this kind can be found in [1--3, 8, 13,
17--19]; see [13] for general references. Moreover, $P$ carries many
algebraic pathologies which also occur elsewhere in abelian group
theory but which can be studied more explicitly in $P$ because of the
sequence-representation of elements in $P$. Such properties quite
often can be phrased in the language of the endomorphism ring of
certain subgroups of products, cf. [5, 7, 9--12, 16, 20].

The most useful tool for constructing such groups is a fundamental
result on $P$ which characterizes slender groups, which is due to
Nunke and based on Specker [21], see [15, Vol.~2, pp.~158--163].  Let
$S=\Bbb Z^{(\aleph_0)}= \oplus_{i \in \omega} e_i \Bbb Z $ denote the
free group on $\aleph_0$ generators, i.e., the canonical direct sum in
the product $P = \prod_{i \in \omega} e_i \Bbb Z$. Recall that a group
$G$ is {\sl slender\/} if any homomorphism from $P$ into $G$ maps
almost all elements $e_i (i \in \omega )$ onto 0. Specker [21] showed
that $\Bbb Z$ is slender. This has many consequence for dual groups,
including the following, where we use the notation $G^*$ for the dual
$\text{Hom}(G,\Bbb Z)$ of $G$.  Specker \cite{21} showed that
$P^*\cong S$, and obviously $S^*\cong P$.  It follows that $P$ has no
direct summand isomorphic to $S$.  More generally, any group $G$ whose
dual is isomorphic to $S$ cannot have a summand isomorphic to $S$, for
$G^*\cong S$ has no summand isomorphic to $S^*\cong P$.  It is natural
to ask how generally ``summand $\cong S$'' is the only obstruction to
``dual $\cong S$.'' Specifically, John Irwin asked whether a pure
subgroup of $P$ that does not have a direct summand isomorphic to $S$
must have a dual isomorphic to $S.$

We shall answer this question negatively, assuming the continuum
hypothesis CH.

The first idea is to phrase the summand question in terms of the
endomorphism ring. Disregarding the condition about the dual, the
existence of subgroups of $P$ with an appropriate endomorphism ring
was established in several papers.  Recall that pure subgroups $G$ of
products of $\Bbb Z$ will always have certain inessential
endomorphisms, which are those with finite rank image plus scalar
multiplication by an integer. The endomorphisms of finite rank
constitute a two-sided ideal $ \text{Fin}\,G$ of $\text{End} \,G$,
where $\text{End}\,G$ denotes the endomorphism ring of $G$. If $G$ can
be chosen so that sums of these and scalars are all the endomorphisms,
i.e., if $\text{End}\,G = \Bbb Z \oplus \text{Fin}\,G$, then the
summand problem is settled because a summand isomorphic to $S$ would
produce many endomorphisms that are not ``almost scalar'' in this
sense.  Such groups are constructed in [6, 9] and inside $P$ in [12,
16, 20]. The new difficulty however is to ensure the condition on the
dual of $G$.  While an easy form of a ``black-box-type argument" as in
[6, 9, 12] will force endomorphisms to be almost scalar, on the other
hand we must ensure that uncountably many homomorphisms from $S$ into
$\Bbb Z$ survive while building $G$ up by a transfinite chain of
extensions from $S$ inside $P$. In order to use topological arguments,
we will restrict our source of new elements for a potential group $G$
to the $\Bbb Z$-adic closure $D$ of $S$ in $P$, which is the pre-image
of the maximal divisible subgroup $D/S$ in $P/S$. It will be important
that this group $D$ can be used to test slenderness, as observed many
years ago by Ti Yen, cf. Fuchs [15, Vol.~2, p.~163, Exercise~5] and
extended recently in [20, p.~276, Corollary~2.5]:

\noindent
$(*)$\quad{\sl A group $X$ is slender if and only if $\text{Hom}(P,X)
= \text{Fin}(P,X)$ if and only if $\text{Hom}(D,X) = \text{Fin}(D,X).$
In particular, every homomorphism from $D$ to $S$ has finite rank.}

Choosing suitable generators by induction we will construct the desired group
$G$ as in the theorem. CH is needed to remain in the ``countable case" at
induction steps.

\proclaim { Theorem}
  Assuming CH, there is a group $G$ with the properties:

 (a) $S \subset  G \subset  D$  and $G$ is pure in $D$ (hence also in
$P).$

(b) $ \text{End} \,G = \Bbb Z \oplus \text{Fin}\,G$, i.e., every
endomorphism of $G$ is almost scalar.

(c)  The dual group $G^* = \text{Hom}(G,\Bbb Z)$ is uncountable.
\endproclaim

\vskip48pt plus16pt minus16pt

\head
{\S 2 Preliminary Facts}
\endhead
We shall need two easy general lemmas.  To avoid interrupting the
proof of the main theorem, we present the lemmas in this preliminary
section.  

We will say that a family $\Cal F$ of subspaces of a vector
space $X$ is {\sl almost independent\/} if it is empty or it
consists of a single subspace different from $X$ or $\Cal F$ is minimal
in the sense that for all  $A \in \Cal F$ we have
 $\sum (\Cal F \setminus \{A\}) \neq \sum \Cal F $.

\proclaim{Lemma 1}
Let $X$ be a rational vector space of infinite dimension $d$, and let
$\Cal F$ be an almost independent family of subspaces of $X$.
Then there are $d$ elements of $X$ that are in none of the spaces in
$\Cal F$.
\endproclaim

Notice that almost independence is important here, even if we have an
upper bound on the cardinality of $\Cal F$ (as we shall, when the
lemma is applied).  $X$ can be covered by countably many proper
subspaces. It is interesting to observe that [7] also needs a lemma
about free space in vector spaces in order to study aspects of $P$.

\demo{Proof of Lemma}
If  $\Cal F = \emptyset $ we have nothing to show and the case
$\Cal F = \{ A \}$ is easy:
If $|A|<d$ then $X \setminus A$ will be the required set,
 by elementary cardinal
arithmetic, since $|X|=d$.  If $|A|=d$, the result is again clear
since $X \setminus A$ contains cosets of $A$, which have cardinality $d$.

Henceforth, suppose $\Cal F$ has at least two members.  Let $A$ be one
of them, and let $B$ be the sum of all the rest.  Both are proper
subspaces of $X$ as $\Cal F$ is almost independent.  If both $A$ and
$B$ had cardinalities $<d$, the result would follow by cardinal
arithmetic.  So assume that one of them, say $A$, has cardinality $d$.
As $A\not\subseteq B$ by almost independence, the argument above (with
$A$ and $\{A\cap B\}$ in place of $X$ and $\Cal F$) shows that there
are $d$ elements $a$ in $A\setminus B$.  Almost independence also
provides an element $b\in B\setminus A$.  The sums $a+b$ for these $d$
values of $a$ and this one value of $b$ are $d$ vectors that are not
in $A\cup B$ and therefore not in $\bigcup\Cal F$.
\qed\enddemo

A {\sl torsionless\/} (abelian) group is (up to isomorphism) a
subgroup of a product $\Bbb Z^\kappa $ or equivalently any group $G$
where $G^*$ separates points in the sense that for any $0 \neq g \in
G$ there is a $\sigma \in G^*$ such that $g \sigma \neq 0$; see [14]
and [15] for more details about torsionless groups.

For any integer $z$ and any group $G$, we write $z$ also for the
endomorphism of $G$ given by multiplication by $z$, and we call such
an endomorphism {\sl scalar\/}. We also write $\Bbb Z$ for the set of
scalar endomorphisms of $G$.  We also use this notation and
terminology when $z$ is a rational number and $G$ is torsion-free, but
then the homomorphism $z$ maps $G$ into its divisible hull rather than
into itself.  We call an endomorphism of $G$ {\sl almost scalar\/} if
it differs from a scalar by a homomorphism which has an image of
finite rank.  The following lemma shows that, when $G$ is torsionless,
it doesn't matter whether we take the scalars to be integers or
rationals in this definition.  (The lemma holds for many groups other
than torsionless ones, but we won't need any more generality.)

\proclaim{Lemma 2}
Let $G$ be a torsionless abelian group of infinite rank.  Let $z$ be a
rational scalar and $\sigma$ a finite-rank homomorphism of $G$ into its
divisible hull.  Suppose $z+\sigma$ maps $G$ into itself.  Then $z$ is an
integer.
\endproclaim

\demo{Proof}
Because $G$ is torsionless and has infinite rank, there is a non-zero
homomorphism $\varphi:G\to\Bbb Z$ vanishing on
$G\cap\text{Range}(\sigma)$.  Without loss of generality, $\varphi$
maps $G$ onto $\Bbb Z$ (just divide by a generator of the range); fix
$e\in G$ with $(e)\varphi=1$.  We write $\varphi$ also for the
homomorphic extension of $\varphi$ mapping the divisible hull of G
into $\Bbb Q$, which vanishes on $\text{Range}(\sigma)$.  Then, as
$e(z+\sigma)\in G$, we have 
$$ 
z=(e)\varphi\cdot
z=(ez)\varphi=(ez+e\sigma)\varphi\in G\varphi=\Bbb Z.  
$$
\qed\enddemo
\bigskip
\head
{\S 3 The Theorem} 
\endhead
\proclaim{Theorem}
Assume the continuum hypothesis.  There is a group $G$ with the
properties:

(a) $S\subseteq G\subseteq D$ and $G$ is pure in $D$ (hence also in $P$).

(b) Every endomorphism of $G$ is almost scalar.

(c) The dual group $G^*=\text{Hom}(G,\Bbb Z)$ is uncountable.
\endproclaim

Before turning to the proof, we point out that the theorem answers
Irwin's question cited in the introduction.  Part (b) of the theorem
implies that $G$ does not have a direct summand isomorphic to $S$, for
such a summand, and therefore $G$ itself, would have many
endomorphisms that are not almost scalar.  And of course part (c)
implies that $G^*$ is not isomorphic to the countable group $S$.

\demo{Proof of Theorem}
$G$ will be obtained as the union of an increasing $\omega_1$-sequence
of subgroups $G_\alpha$ ($\alpha<\omega_1$), where $G_0=S$,
$G_\lambda=\bigcup_{\alpha<\lambda}G_\alpha$ for limit $\lambda$, and
each $G_\alpha$ is a countable, pure subgroup of $D$.  This ensures
that part (a) of the theorem holds.  The successor stages of the
induction are designed to ensure parts (b) and (c).

To deal with (b), we begin with an enumeration, in an
$\omega_1$-sequence $(\eta_\alpha)_{\alpha<\omega_1}$, of all
homomorphisms $\eta:S\to P$.  We require further that each such
homomorphism $\eta$ occurs $\aleph_1$ times in the enumeration.  The
continuum hypothesis gives us such an enumeration, because there are
$2^{\aleph_0}$ such homomorphisms.  Each $\eta_\alpha$ has a unique largest
extension with domain pure in $D$, say $\be:D_\alpha\to P$.  This is
obtained by first extending $\eta_\alpha$ to a homomorphism from $D$,
which is the $\Bbb Z$-adic closure of $S$ in $P$, to the $\Bbb Z$-adic
completion of $P$, and then restricting to the subgroup $D_\alpha$ of
$D$ that is mapped into $P$.

At stage $\alpha$ of the construction, when we have $G_\alpha$ and are
defining $G_{\alpha+1}$, we shall ensure that $\be$ does not map
$G_{\alpha+1}$ into $G$ unless $\be\res G_\alpha$ is almost scalar.
If we achieve this, then (b) will hold.  To see this, suppose $e$ were
an endomorphism of $G$ that is not almost scalar.  For each scalar
$z$, the endomorphism $e-z$ of $G$ has infinite rank, so there is an
ordinal $\beta<\omega_1$ so large that $e-z\res G_\beta$ has infinite
rank.  As there are only countably many scalars, we can fix one
$\beta$ that works for all $z$.  Now $e\res S:S\to P$ occurs
uncountably often in the enumeration $(\eta_\alpha)$, so fix an
$\alpha\geq\beta$ for which $e=\eta_\alpha$.  Then, at stage $\alpha$
of the construction, since $\be\res G_\alpha=e\res G_\alpha$, we
ensure that $\be$ does not map $G_{\alpha+1}$ into $G$.  This is
absurd, since $\be$, the largest extension of $\eta_\alpha$ to a pure
subgroup of $D$, extends $e$, which maps all of $G$ into $G$.

If $\be\res G_\alpha$ is almost scalar or if the domain of $\be$ does
not include all of $G_\alpha$, then we need not do anything at stage
$\alpha$ and we set $G_{\alpha+1}=G_\alpha$.  Otherwise, we shall set
$G_{\alpha+1}=\langle G_\alpha\cup\{x\}\rangle_*$, the pure subgroup
of $D$ generated by $G_\alpha$ and one new, carefully chosen $x\in D$.
If $x$ is in the domain $D_\alpha$ of $\be$, then we shall also ensure
that $x\be\notin G_\beta$ for all $\beta$; thus $\be$ will not map
$G_{\alpha+1}$ into $G$.  Of course, in order to do this, we must
choose $x$ so that $x\be\notin G_{\alpha+1}$ and, at later stages of
the construction, we must be careful not to put $x\be$ into any
$G_\beta$.

Thus, as the construction proceeds, we accumulate a set $V$ of
forbidden elements, which must never be put into any $G_\beta$.  Each
stage contributes at most one element $x\be$ to $V$, so $V$ is
countable at each stage of the construction.  Of course the choice of
$x$ at stage $\alpha$ is constrained by the elements already in $V$
from previous stages.

To deal with part (c) of the theorem, we shall have, at each stage
$\alpha$, some homomorphisms $f^\alpha_\xi:G_\alpha\to\Bbb Z$, indexed
by the ordinals $\xi<\alpha$, such that for $\xi<\alpha<\beta$ we have
$f^\alpha_\xi=f^\beta_\xi\res G_\alpha$.  Then, for each
$\xi<\omega_1$, the homomorphisms $f^\alpha_\xi$ for all $\alpha$
(from $\xi+1$ up) combine to give a homomorphism $f_\xi:G\to\Bbb Z$.
If we make the $f^\alpha_\xi$ at each stage $\alpha$ distinct (for
distinct $\xi$), then all these $f_\xi$ will be distinct and (c) will
hold.

Thus, what needs to be done at stage $\alpha$ of the construction is
the following.  We are given

\roster
\item a countable pure subgroup $G_\alpha$ of $D$,
\item countably many homomorphisms $f^\alpha_\xi:G_\alpha\to\Bbb Z$,
indexed by $\xi<\alpha$,
\item a countable set $V_\alpha$ of forbidden elements, with
$V_\alpha\cap G_\alpha=\emptyset$, and
\item a homomorphism $\be:D_\alpha\to P$ whose domain includes
$G_\alpha$ and which is not almost scalar on $G_\alpha$.
\endroster
We seek an element $x\in D$ such that
\roster
\item[5] $\langle G_\alpha\cup\{x\}\rangle_*$ contains neither
$x\be$ nor any element of $V_\alpha$, and
\item each $f^\alpha_\xi$ extends to $f^{\alpha+1}_\xi:\langle
G_\alpha\cup\{x\}\rangle_*\to\Bbb Z$.
\endroster
In addition, we need one new homomorphism $f^{\alpha+1}_\alpha:\langle
G_\alpha\cup\{x\}\rangle_*\to\Bbb Z$, distinct from all
$f^{\alpha+1}_\xi$ for $\xi<\alpha$.

Getting $f^{\alpha+1}_\alpha$ is easy, as $\langle
G_\alpha\cup\{x\}\rangle_*$ is a countable subgroup of $P$ containing
$S$, hence is free of infinite rank by $(\dag)$ of the introduction,
and hence has uncountably many homomorphisms to $\Bbb Z$.  So our
efforts from now on will be directed toward finding an $x\in D$
subject to \therosteritem5 and
\therosteritem6.  Once we have such an $x$, we can set
$G_{\alpha+1}=\langle G_\alpha\cup\{x\}\rangle_*$ and
$V_{\alpha+1}=V_\alpha\cup\{x\be\}$, and the proof will be complete.

Notice that $G_\alpha$, being a countable, pure subgroup of $P$ and
dense in the product topology of $P$, can be mapped onto $S$ by an
automorphism of $P$, according to a theorem of Chase [4, p.~605,
Corollary~3.3]. Furthermore, since $P/S=(D/S)\oplus R$ where $R$ is
reduced and $D/S$ is divisible, and since $S\subseteq
G_\alpha\subseteq D$, we have $P/G_\alpha=(D/G_\alpha)\oplus R$ with
$D/G_\alpha$ divisible (being a quotient of $D/S$).  So $D$ is the
pre-image in $P$ of the divisible part of $P/G_\alpha$, as well as the
pre-image in $P$ of the divisible part of $P/S$.  Therefore, an
automorphism of $P$ that maps $G_\alpha$ onto $S$ necessarily maps $D$
onto $D$.  Applying such an automorphism, we obtain the following
description of what needs to be done at stage $\alpha$ (numbered to
match the earlier description).  We are given
\roster
\item[2] countably many homomorphisms $f_i:S\to\Bbb Z$
($i \in\omega$), where we have re-indexed by $\omega$ in place of
$\alpha$,
\item a countable set $V\subseteq P-S$,
\item a homomorphism $\eb:E\to P$, whose domain $E$ includes $S$, such
that $\eb\res S$ is not almost scalar.
\endroster
We seek an $x\in D$ such that
\roster
\item[5] $\langle S\cup\{x\}\rangle_*$ contains neither $x\eb$ nor
any member of $V$, and
\item each $f_i$ extends to a homomorphism $\langle
S\cup\{x\}\rangle_*\to\Bbb Z$.
\endroster

The $V$ part of \therosteritem5 is easy to handle.  We must exclude
any $x$ such that
$$
s+xm=vn\tag7
$$
for some $s\in S$, $m\in\Bbb Z$, $n\in\Bbb Z-\{0\}$, and $v\in V$.  We
need not worry about $m=0$, as $S\cap V=\emptyset$ (and $S$ is pure).
Then, for each $s$, $m$, $n$, and $v$, there is at most one $x$
satisfying \thetag7 (as $P$ is torsion-free).  So only countably many
$x$ are excluded.  Therefore, it suffices to find uncountably many
$x\in D$ satisfying \therosteritem6 and
$$
x\eb\notin\langle S\cup\{x\}\rangle_*.\tag5'
$$
Here \thetag{$5'$} is regarded as vacuously true if $x\notin E$.

Fix an enumeration $(r_n)_{n\in\omega}$ of the rational numbers, such
that each rational number occurs infinitely often in the enumeration.
We shall construct a sequence of elements $b_n\in S-\{0\}$ with the
following properties.  We write $\supp(x) = \{ i \in \omega : x_i \neq
0 \}$ for the set of indices where an element $x = \sum_{i \in \omega}
e_i x_i\in P$ has non-zero components $x_i$; for elements of $S$, this
is a finite set, so the maximum and minimum mentioned in
\therosteritem8 below make sense.
\roster
\item[8] If $m<n$ then $\max\supp(b_m)<\min\supp(b_n)$.
\item If $i\leq n$ then $b_n f_i=0$.
\item $b_n \eb - b_n r_n$ has at least one non-zero component beyond
position $n$.
\endroster
In \therosteritem{10}, we multiplied by a rational scalar, so we
should work in the divisible hull of $P$.  Alternatively, we can
reformulate such statements (here and below) to involve only integer
multipliers, by clearing the denominators.

The $b_n$ are defined inductively.  Suppose $b_0,\dots,b_{n-1}$ are
given, and let $q\in\omega$ be larger than all elements of their
supports.  We shall choose $b_n$
 in the subgroup $S_q = \bigoplus_{i \geq q} e_i \Bbb Z$
 of $S$
consisting of elements of $S$ whose supports have minimum $\geq q$.
That will ensure that \therosteritem8 holds.  Since $S_q$ has finite
corank $q$ in $S$, the range of $(\eb-r_n)\res S_q$ has infinite rank.
(Otherwise, $\eb\res S$ would be the scalar $r_n$ plus an endomorphism
of finite rank, contrary to \therosteritem4.)  So we can find a
finitely generated subgroup $F$ of $S_q$ whose image under $\eb-r_n$
has rank $\geq 2n+3$.  The subgroup $F'$ of $F$ consisting of elements
$b$ with $(\forall i\leq n)\,b f_i=0$, i.e., the subgroup of elements
that could serve as $b_n$ and satisfy \therosteritem9, has corank
$\leq n+1$ in $F$ (being the intersection of the kernels of $n+1$
homomorphisms $f_i$ to $\Bbb Z$).  So the image of $F'$ under
$\eb-r_n$ has rank at least $(2n+3)-(n+1)>n+1$.  This image must
therefore contain an element with a non-zero component past position
$n$.  That is, we can choose $b_n\in F'$ (thereby satisfying
\therosteritem9) in such a way that \therosteritem{10} also holds.
This completes the construction of the sequence $(b_n)$.

Define an endomorphism $\beta$ of $P$ taking each $e_i$ to $b_i$, that is
$$
\beta:P\to P: \sum_{i=0}^\infty e_ix_i \mapsto\sum_{i=0}^\infty b_ix_i.
$$
The infinite sum in this definition makes sense because the $b_i$
have disjoint supports by \therosteritem8.  It is trivial to check
that $\beta$ maps $S$ into $S$ and $D$ into $D$.  (The latter also
follows from the former by $\Bbb Z$-adic continuity.)  Since all the
$b_i$ are non-zero, $\beta$ is one-to-one.

Also define, for each $i\in\omega$, a partial homomorphism
$$
g_i:P\rightharpoonup\Bbb Z: \sum_{j=0}^\infty e_jx_j
 \mapsto\lim_{n\to\infty}
(\sum _{i=0}^{q_n} e_j x_j) f_i ,
$$
where $q_n=\max\supp(b_n)$, and where the limit means the value for
all sufficiently large $n$ provided this value is eventually
independent of $n$.  (If the value is not eventually independent of
$n$, the limit is undefined; that's why $g_i$ is only a partial
homomorphism.)  The domain of $g_i$ is clearly a pure subgroup of $P$
that includes $S$, and $g_i\res S=f_i$.  Furthermore, if
$ \sum_{i=0}^\infty e_ix_i = (\sum_{i=0}^\infty e_iy_i)\beta$, then
$$
(\sum _{i=0}^{q_{n+1}} e_i x_i) f_i -(\sum _{i=0}^{q_n} e_i x_i) f_i
=( \sum _{i=q_n+1}^{q_{n+1}} e_i x_i) f_i
=(b_{n+1}y_{n+1})f_i = 0
$$
by \therosteritem9 provided $n\geq i$.  Thus
$P\beta\subseteq\text{domain}(g_i)$.  We shall ensure \thetag6 by
choosing $x$ to be $y\beta$ for some $y\in D$; then $\langle
S\cup\{x\}\rangle_*\subseteq\text{domain}(g_i)$, so we can extend
$f_i$ to $g_i\res\langle S\cup\{x\}\rangle_*$.

To obtain uncountably many $x\in D$ satisfying \thetag{$5'$} and
\therosteritem6, which is all we need to finish the proof, it now
suffices, as $\beta$ is one-to-one, to find uncountably many $y\in D$
such that
$$
y\beta\eb\notin\langle S\cup\{y\beta\}\rangle_*.\tag{11}
$$

There is an easy case, namely if
$D\beta\not\subseteq\text{domain}(\eb)$, so $\beta\eb$ is not
defined on all of $D$.  Then the domain $E\beta^{-1}$ of
$\beta\eb$ is a pure, proper subgroup of $D$ and includes $S$ (as
$S\beta\subseteq S\subseteq E$).  So $E\beta^{-1}/S$ is a proper
subspace of the uncountable-dimensional rational vector space $D/S$.
By Lemma~1, the set-theoretic complement of such a subspace is
uncountable.  And any $y\in D$ whose image is in that complement
satisfies \thetag{11} vacuously.  This finishes the easy case.

Henceforth, assume that $\beta\eb$ is defined on all of $D$.  To
find uncountably many $y\in D$ satisfying \thetag{11}, we consider how
some $y$ could fail to satisfy \thetag{11}.  That would mean that
$$
y\beta\eb\cdot q= y\beta \cdot p + s
$$
for some $s\in S$ and some integers $p,q$ with $q\neq0$.  Write $r$
for the rational number $p/q$.  Then
$y\beta\eb-y\beta\cdot r = s(1/q)$ has finite support.

For each rational number $r$, let
$$
W_r=\{y\in D\mid
y\beta\eb-y\beta \cdot r \text{ has finite support}\}.
$$
The preceding discussion shows that any $y\in D$ violating \thetag{11}
is in some $W_r$.  So, to complete the proof, we need to find
uncountably many elements of $D$ that are in no $W_r$.

Since both $\eb$ and $\beta$ map $S$ into itself, we clearly have
$S\subseteq W_r$.  Also, $W_r$ is a pure subgroup of $D$.  Therefore,
$W_r/S$ is a vector subspace of the rational vector space $D/S$.  We
shall prove that these subspaces are proper, linearly independent,
 hence almost independent subspaces of $D/S$.
   Once this is done, Lemma~1 will produce
uncountably many elements of $D/S$ not contained in any $W_r$; the
pre-images in $D$ of these elements will be the desired $y$'s, and the
proof will be complete.

First, we show that $W_r/S$ is a proper subspace of $D/S$.  Suppose
not, i.e., suppose $W_r=D$ for a certain $r=p/q$, where $p$ and $q$
are integers.  In view of the definition of $W_r$,
$$
y\mapsto (y\beta\eb-y\beta \cdot r) \cdot q
=y\beta\eb \cdot q -y\beta \cdot p
$$
is a homomorphism from $W_r=D$ into $S$.  But all homomorphisms $D\to
S$ have finite rank by (*) in the introduction.
  So $D\beta(\eb-r)$ would have finite rank.
But this group contains, by definition of $\beta$, the elements
$b_n(\eb-r)$ for all $n$, including the infinitely many $n$ for
which $r_n=r$.  But those elements generate a group of infinite rank,
by property \thetag{10} of the $b_n$.  This contradiction shows that
$W_r/S$ is a proper subspace of $D/S$.

It remains to show that the subspaces $W_r/S$ are linearly
independent.  This proof is a minor variation of the standard proof
that eigenspaces of a linear operator corresponding to different
eigenvalues are linearly independent.  Suppose we had a linear
dependence between these subspaces of $D/S$.  Back in $D$, this would
mean
$$
\sum_{i=1}^k w_i c_i\in S,\tag{12}
$$
where $c_i\in\Bbb Q-\{0\}$, $k\geq1$, $w_i\in W_{r_i}-S$, and all
$r_i$ are distinct.  (As before, we adopt the convention that we
either work with divisible hulls and in particular that $S$ in
\thetag{12} really means its divisible hull, so that the rational
scalars make sense, or interpret statements like \thetag{12} as
meaning the result of clearing denominators, so that only integer
scalars occur.)  Take a relation of the form \thetag{12} with $k$ as
small as possible.  Note that $k\geq2$, for $w_1c_1\in S$ would imply
$w_1\in S$, a contradiction.  Apply $\beta\eb-\beta\cdot r_1$ to
\thetag{12}, remembering that $S$ is closed under $\eb$ and $\beta$
and that $\beta\eb$ is defined on all of $D$ so that this makes sense.
The result is, in view of the definition of $W_{r_i}$,
$$
\sum_{i=1}^k[w_i\beta (r_i-r_1)+\text{ elements of }S] c_i\in S.
$$
Therefore,
$$
\left(\sum_i w_i(r_i-r_1)c_i\right)\beta=\sum_i w_i\beta (r_i-r_1)(c_i)
\in S.
$$
It is immediate from the definition of $\beta$ (since no $b_i$ is 0)
that $S\beta^{-1}=S$.  So we have
$$
\sum_i w_i(r_i-r_1)c_i\in S.
$$
The $i=1$ term here contains $r_1-r_1$, so it vanishes, but the other
terms have non-zero coefficients.  So we have a relation of the form
\thetag{12} with fewer summands.  This contradicts the minimality of
$k$, and this contradiction completes the proof.
\qed\enddemo

\enddocument